\documentclass[11pt]{amsart}
\usepackage{graphicx}
\usepackage{amsmath,color,graphicx}
\usepackage{amsfonts}
\usepackage{amsthm}
\usepackage{amsbsy}
\usepackage{amssymb}
\usepackage{microtype}
\usepackage{amsrefs}
\usepackage{hyperref}
\usepackage{mathtools}
\usepackage{multicol}
\usepackage{cancel}
\usepackage{xcolor}
\usepackage{bookmark}
\usepackage{enumerate}
\usepackage[margin=1in,footskip=0.25in]{geometry} 
\numberwithin{equation}{section}
\newtheorem{thm}{Theorem}[section]
\theoremstyle{definition}
\newtheorem{lem}[thm]{Lemma}
\newtheorem{prop}[thm]{Proposition}

\newtheorem{cor}[thm]{Corollary}

\newtheorem{ques}[thm]{Question}
\newtheorem{defn}[thm]{Definition}

\newtheorem{rem}[thm]{Remark}
\newtheorem{exam}[thm]{Example}
\newtheorem*{exam*}{Example}

\newcommand{\B}{{\mathcal{B}}}

\newcommand{\F}{{\mathcal{F}}}

\renewcommand{\phi}{\varphi}

\usepackage{graphicx}
\usepackage{amsmath,color,graphicx}
\usepackage{amsfonts}
\usepackage{amsthm}
\usepackage{amsbsy}
\usepackage{amssymb}
\usepackage{microtype}
\usepackage{amsrefs}
\usepackage{hyperref}
\usepackage{mathtools}
\usepackage{multicol}
\usepackage{cancel}
\usepackage{xcolor}
\usepackage{bookmark}
\usepackage{tikz}
\usepackage{xcolor}
\definecolor{darkgreen}{RGB}{0,100,0}

\title[Bruhat Shellings of the Quasisymmetric Flag Moment Complex]{Bruhat Shellings of the\\ Quasisymmetric Flag Moment Complex}
\author{Saffana Khan}
\address{Department of Mathematics and Statistics, York University, Toronto, Ontario M3J 1P3, Canada}
\email{\href{mailto:khansaff@my.yorku.ca}{khansaff@my.yorku.ca}}

\author{Nantel Bergeron}
\address{Department of Mathematics and Statistics, York University, Toronto, Ontario M3J 1P3, Canada}
\email{\href{mailto:bergeron@yorku.ca}{bergeron@yorku.ca}}
\urladdr{https://bergeron.mathstats.yorku.ca/}
\date{July 2026}
\begin{document}

\begin{abstract}
The quasisymmetric flag variety gives rise to a regular cell complex whose
maximal cells are cubes indexed by fully supported noncrossing partitions.
Two maximal cubes can meet along several faces.
We describe each cube using the block arcs and parent connectors of
its indexing partition, and determine its intersection with the union of
the preceding cubes in Bruhat order.
For each nonempty face, we determine the first maximal cube containing it by
applying an outer-block merging operation to its Bruhat-maximal
noncrossing-partition label. The first-carrier theorem provides an explicit description of how each cube attaches to the preceding subcomplex. More precisely, the intersection is a union of coordinate facets which, for every cube except the first, forms the closed star of a face in its boundary. Consequently, every linear extension of Bruhat order on the maximal-cell labels gives a shelling of the regular cell complex. In particular, its barycentric subdivision is shellable, hence homotopy Cohen--Macaulay and Cohen--Macaulay over \(\mathbb Z\) and every field.
\end{abstract}
\maketitle

\section{Introduction}

Borel identified the integral cohomology of the variety
\(\operatorname{Fl}_n\) of complete flags in \(\mathbb C^n\) with the
coinvariant ring of the symmetric group~\cite{Borel1953}. Replacing
the symmetric relations by all quasisymmetric polynomials with zero
constant term gives the quasisymmetric coinvariant ring. Over \(\mathbb Q\)
its dimension is the \(n\)th Catalan number, the number of noncrossing
partitions of \(\{1,\ldots,n\}\)~\cite{AvalBergeronBergeron2004}.
Bergeron, Gagnon, Nadeau, Spink, and Tewari realized the integral quasisymmetric coinvariant ring as the
integral cohomology of a subvariety
\(\operatorname{QFl}_n\subseteq\operatorname{Fl}_n\), the
\emph{quasisymmetric flag variety}. Its fixed points under the diagonal
torus are indexed by those same noncrossing partitions~\cite[Theorem~A and
Section~6]{BergeronGagnonNadeauSpinkTewari2026}.
Each irreducible component is toric for this action, and its moment polytope
has the face structure of a cube, with one nonempty face for each torus
orbit closure in the component. Identifying faces of different cubes
whenever they represent the same orbit closure gives the \emph{moment
complex} \(\mathcal C_n\)~\cite[Section~8.3]{BergeronGagnonNadeauSpinkTewari2026}.

We study how the maximal cubes of \(\mathcal C_n\) fit together. Even two squares in this complex can share two edges. As we add the cubes one at a time, which part of each new cube is already present? For any ordering of the cubes that respects Bruhat order on their labels, we show how to read this intersection directly from the noncrossing partition indexing the new cube. This gives a shelling of \(\mathcal C_n\) and proves that its barycentric subdivision is Cohen--Macaulay.

For \(n=3\) there are just two maximal squares, indexed by the noncrossing partitions \(123\) and \(13|2\) of \(\{1,2,3\}\), where a bar separates blocks. In each diagram below, black arcs join successive elements of the same block. In \(13|2\), the singleton \(2\) lies under the arc from \(1\) to \(3\). The dashed blue arc above joins \(2\) to the right endpoint \(3\) of that arc. We call it a \emph{connector}.

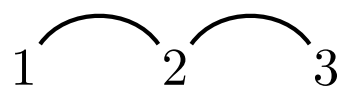
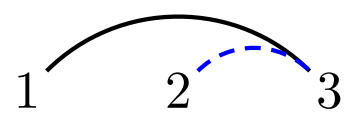
\[
\begin{tikzpicture}[x=1cm,y=1cm,color=black,font=\small]
    \begin{scope}
        \foreach \i in {1,2,3}
            \node[inner sep=1pt] (a\i) at (\i,0) {\(\i\)};
        \draw[thick] (a1) to[out=55,in=125] (a2);
        \draw[thick] (a2) to[out=55,in=125] (a3);
        \node at (2,-0.72) {\(123\)};
    \end{scope}
    \begin{scope}[xshift=4.5cm]
        \foreach \i in {1,2,3}
            \node[inner sep=1pt] (b\i) at (\i,0) {\(\i\)};
        \draw[thick] (b1) to[out=45,in=135] (b3);
        \draw[blue,densely dashed,thick] (b2) to[out=45,in=135] (b3);
        \node at (2,-0.72) {\(13|2\)};
    \end{scope}
\end{tikzpicture}
\]

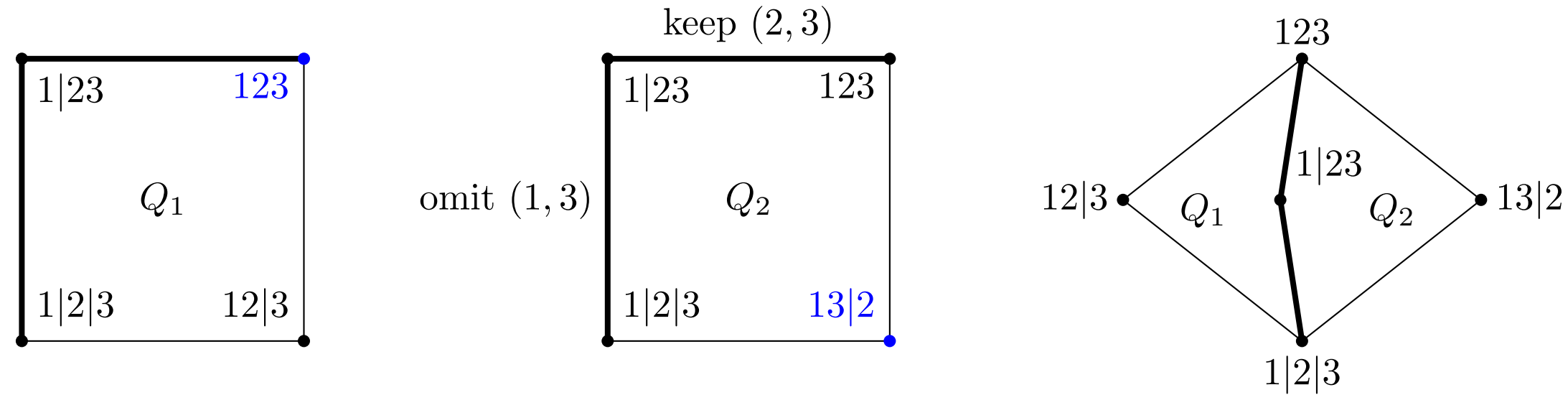
\begin{figure}[htbp]
\centering
\begin{tikzpicture}[x=1cm,y=1cm,font=\small]
    \colorlet{ordinaryvertex}{.}
    \foreach \shiftx/\cubeindex/\cornerlabel/\uppercolor/\lowercolor in
        {0/1/{12|3}/blue/ordinaryvertex,5.4/2/{13|2}/ordinaryvertex/blue} {
        \begin{scope}[xshift=\shiftx cm]
            \draw (0,0) rectangle (2.6,2.6);
            \draw[line width=1.5pt] (0,0) -- (0,2.6) -- (2.6,2.6);
            \foreach \a/\b in {0/0,0/2.6}
                \fill (\a,\b) circle (1.6pt);
            \node[above right] at (0,0) {\(1|2|3\)};
            \node[below right] at (0,2.6) {\(1|23\)};
            \fill[\uppercolor] (2.6,2.6) circle (1.6pt);
            \node[below left,text=\uppercolor] at (2.6,2.6) {\(123\)};
            \fill[\lowercolor] (2.6,0) circle (1.6pt);
            \node[above left,text=\lowercolor] at (2.6,0) {\(\cornerlabel\)};
            \node at (1.3,1.3) {\(Q_{\cubeindex}\)};
        \end{scope}
    }
    \node[left] at (5.4,1.3) {omit \((1,3)\)};
    \node[above] at (6.7,2.6) {keep \((2,3)\)};

    \begin{scope}[xshift=11.8cm,line join=round]
        \coordinate (bottom) at (0,0);
        \coordinate (middle) at (-0.20,1.3);
        \coordinate (top)    at (0,2.6);
        \coordinate (west)   at (-1.65,1.3);
        \coordinate (east)   at (1.65,1.3);
        \draw (bottom) -- (west) -- (top) -- (east) -- cycle;
        \draw[line width=1.5pt] (bottom) -- (middle) -- (top);
        \foreach \v in {bottom,middle,top,west,east}
            \fill (\v) circle (1.6pt);
        \node[below]       at (bottom) {\(1|2|3\)};
        \node[above right] at (middle) {\(1|23\)};
        \node[above]       at (top)    {\(123\)};
        \node[left]        at (west)   {\(12|3\)};
        \node[right]       at (east)   {\(13|2\)};
        \node at (-0.91,1.20) {\(Q_1\)};
        \node at ( 0.83,1.20) {\(Q_2\)};
    \end{scope}
\end{tikzpicture}
\caption{The two maximal squares of \(\mathcal C_3\).
Bars separate blocks. Each cube is indexed by its unique maximum vertex label in Bruhat order. The first two panels show these maxima in blue, \(123\) for \(Q_1\) and \(13|2\) for \(Q_2\). The third panel shows the squares glued along the bold path from \(1|2|3\) through \(1|23\) to \(123\).}
\label{fig:n3-common-path}
\end{figure}

To get a vertex of the corresponding square, keep or omit each of the two edges independently. The connected components give the blocks of the vertex label, and changing one choice moves along an edge of the square. In the diagram for \(13|2\), keeping neither edge gives \(1|2|3\), keeping only the connector gives \(1|23\), and keeping both gives \(123\). These three vertices also occur in the first square. Keeping only the black arc gives the one new vertex, \(13|2\). The two edges of the second square that avoid this vertex form the common path shown in bold in Figure~\ref{fig:n3-common-path}.

The same construction works for any noncrossing partition \(P\) of \([n]=\{1,\ldots,n\}\). Each arc or connector gives one independent choice in a cube \(Q_P\), and faces with the same set of vertex labels are identified. The maximal cubes have dimension \(n-1\). Their labels are the \emph{fully supported} partitions, those in which \(1\) and \(n\) belong to the same block. We order the labels by writing each block as a decreasing cycle and comparing the resulting permutations in Bruhat order. For example, \(123\) and \(13|2\) give the permutations \(312\) and \(321\) in one-line notation. Since \(312<_{B}321\) in Bruhat order, the square indexed by \(123\) comes first.

Now list the fully supported labels \(\mu_1,\ldots,\mu_N\) in any linear extension of Bruhat order and set
\[
    Q_i:=Q_{\mu_i},\qquad
    Y_{j}:=\bigcup_{1\le i\le j} Q_{i},
    \qquad Y_0:=\varnothing.
\]

\begin{thm}[Shellability]\label{thm:shelling-intro}
Every linear extension of Bruhat order on the top-cell labels is a shelling order of \(\mathcal C_n\).
\end{thm}

For a regular cell complex \(K\), its barycentric subdivision
\(\operatorname{sd}(K)\) is the simplicial complex whose vertices are the
nonempty faces of \(K\) and whose simplices are chains of such faces.

\begin{cor}[Homotopy Cohen--Macaulayness]\label{cor:cm-intro}\label{cor:homotopy-cm}
For every \(1\le j\le N\), the barycentric subdivision
\(\operatorname{sd}(Y_j)\) is shellable and hence
homotopy Cohen--Macaulay. In particular, it is Cohen--Macaulay over
\(\mathbb Z\) and every field.
\end{cor}

The conclusion controls local topology: every link of dimension
\(r\) in \(\operatorname{sd}(Y_j)\) is \((r-1)\)-connected.
This remains true despite branching among the cubes. Already for \(n=4\),
three maximal cubes share a square, and \(\mathcal C_4\) is not a manifold
with boundary (Example~\ref{ex:three-carriers}).

The shelling follows from two descriptions of the way faces enter the
complex. We call a maximal cube containing a face \(H\) a
\emph{carrier} of \(H\). Every nonempty face \(H\) has a largest vertex
label \(\lambda(H)\) in Bruhat order
(Proposition~\ref{prop:partition-model}).
Given a partition \(P\), let \(\kappa(P)\) be obtained by merging
its outer blocks, those that do not lie under an arc of another block.
This operation determines the first maximal carrier.

\begin{thm}[First-carrier theorem]\label{thm:first-carrier-intro}
Let \(H\) be a nonempty face of \(\mathcal C_n\) and put \(P=\lambda(H)\). Among the maximal cubes containing \(H\), the one with least label in Bruhat order is
\(
    Q_{\kappa(P)}.
\)
It is therefore the first carrier of \(H\) in every linear extension of Bruhat order on the maximal-cube labels.
\end{thm}

Merging the outer blocks adds an arc between each consecutive pair
of outer blocks. Omitting these new arcs recovers \(Q_P\) as a face of
\(Q_{\kappa(P)}\). The first-carrier description then determines the
intersection with the preceding union.

\begin{thm}[Attaching-locus theorem]\label{thm:attaching-intro}
For \(j>1\), there is a face \(R\) of \(Q_j\) whose vertices are precisely the new vertices at step \(j\). The intersection \(Q_j\cap Y_{j-1}\) is the union of precisely those facets of \(Q_j\) disjoint from \(R\).
\end{thm}

In Figure~\ref{fig:n3-common-path}, the face \(R\) is the single
new vertex \(13|2\). For \(\mu=134|2\), it is an edge: keep the arc
\((1,3)\), omit the connector \((2,3)\), and vary the choice of the
outer unit arc \((3,4)\). Its endpoints are \(13|2|4\) and \(134|2\).

In the square indexed by \(13|2\), one shared edge is obtained by omitting the black arc over \(2\) and the other by keeping the blue connector. In every dimension, keeping any connector gives an attaching facet. Omitting a black arc gives one precisely when it passes over a vertex or belongs to a nested block (Theorem~\ref{thm:attaching}). Together these facets form the closed star of the face opposite \(R\) in the cube boundary. They can therefore come first in a shelling of that boundary.

An earlier non-pure Bruhat paving is indexed by all noncrossing partitions
\cite[Theorem~9.1, Remark~9.2, and Lemma~9.3]{BergeronGagnonNadeauSpinkTewari2026}.
In the Coxeter flag construction, the corresponding moment polytopes
attach along boundary faces avoiding a distinguished simple vertex;
this also proves contractibility
\cite[Theorem~9.3 and Remark~9.17]{CoxeterFlag2026}.
Those polytopes have varying dimensions. Here we group the pieces by their
least maximal carrier and describe the attaching loci of the maximal cubes.
Their new vertices form a face that may have positive dimension, and
the resulting order is a pure shelling.

Bruhat order on noncrossing partitions and its relation to full support
were originally studied by Gobet--Williams \cite{GobetWilliams2016} and
Biane--Josuat-Verg\`es \cite{BianeJosuatVerges2019}.

Appendix~\ref{app:comparison} relates the partition diagrams back to
\(\operatorname{QFl}_n\). It identifies the vertex labels of each face
with the torus fixed points on the corresponding orbit closure.
Thus the shelling described above applies to the geometric moment complex. We close by asking whether the local rings of the variety \(\operatorname{QFl}_n\) are Cohen--Macaulay as well (Question~\ref{ques:geometric-cm}).

\section{Noncrossing partitions and their cubes}

\begin{defn}[Noncrossing partitions]
A noncrossing partition of \([n]=\{1,\ldots,n\}\) is a set partition
    \(P=\{A_1,\ldots,A_k\}\) with \(A_1\sqcup\cdots\sqcup A_k=[n]\) such that, for any two distinct blocks \(A_i\) and \(A_j\), there are no elements \( a,b\in A_i\) and \(c,d\in A_j\) with \(a<c<b<d.\)
We write \(\operatorname{NC}_n\) for the set of these partitions.

We represent a block \(B=\{b_1<\cdots<b_m\}\) by its consecutive arcs
\[
    E_{\mathrm{arc}}(B)
      :=\{(b_1,b_2),\ldots,(b_{m-1},b_m)\}.
\]
Following the conventions of \cite{BergeronGagnonNadeauSpinkTewari2025,BergeronGagnonNadeauSpinkTewari2026}, we use the standard identification between noncrossing set partitions and noncrossing permutations, under which a block \(B=\{b_1<\cdots<b_m\}\) corresponds to the decreasing cycle \((b_m\,\cdots\,b_1)\). Thus, we use the same symbol \(P\) for a noncrossing set partition and its associated noncrossing permutation when no confusion can arise. Throughout, \(\le_B\) denotes Bruhat order on these permutations.
\end{defn}

\begin{exam}\label{ex:nc-partition-permutation}
The partition
\(
    P=\bigl\{\{1,4,5\},\{2,3\},\{6\},\{7,8\}\bigr\}
\)
is noncrossing. In the notation of Figure~\ref{fig:n3-common-path}, we write
\(P=145|23|6|78\). We omit braces and commas, write the elements of each
block in increasing order, and separate the blocks by vertical bars. For multi-digit entries, we retain commas within blocks whenever needed to avoid ambiguity.
Its consecutive arcs are \((1,4)\), \((4,5)\),
\((2,3)\), and \((7,8)\). Reading each block as a decreasing cycle gives
\(\textcolor{blue}{P=(5\,4\,1)(3\,2)(6)(8\,7)}\).
The diagram shows each value \(P(i)\) in blue below the corresponding label \(i\).
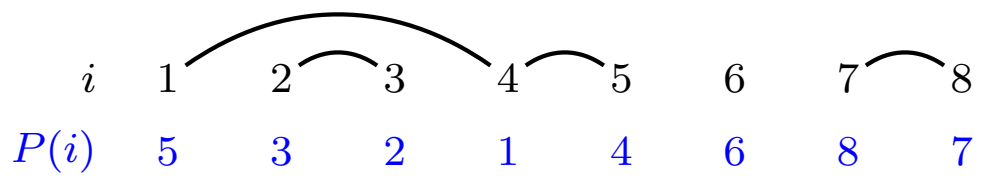
\[
\begin{tikzpicture}[scale=0.75]
    \foreach \i in {1,...,8} {
        \node[inner sep=1pt] (v\i) at (\i-1,0) {$\scriptstyle\i$};
    }
    \foreach \i/\j in {1/4,4/5,2/3,7/8} {
        \draw[thick] (v\i) to[out=35,in=145] (v\j);
    }
    \node[anchor=east] at (-0.45,0) {$\scriptstyle i$};
    \node[blue,anchor=east] at (-0.45,-0.65) {$\scriptstyle P(i)$};
    \foreach \i/\pi in {1/5,2/3,3/2,4/1,5/4,6/6,7/8,8/7} {
        \node[blue] at (\i-1,-0.65) {$\scriptstyle\pi$};
    }
\end{tikzpicture}
\]
\end{exam}

\textbf{Nesting terminology.}
Let \(A,B\) be distinct blocks of \(P\in\operatorname{NC}_n\). We say that \(B\) is \emph{nested inside} \(A\) if
\[
    \min A<\min B\le \max B<\max A.
\]
A block is \emph{outer} if it is not nested inside any other block. Further, we call \(P\) \emph{fully supported} if it has exactly one outer block. Equivalently, \(1\) and \(n\) belong to the same block.

\begin{defn}\label{def:kappa}
    Let \(P\in\operatorname{NC}_n\), and let \(B_1, \dots, B_t\) denote its outer blocks, indexed so that
    \(
        \max B_i < \min B_{i+1},
    \) for  $1\le i<t$.
    Define \(\kappa(P)\) to be the noncrossing set partition obtained by merging all outer blocks of \(P \) into a single block, while leaving every nonouter block unchanged. Equivalently, 
    \[
        \kappa(P) = \{B_1 \cup \cdots \cup B_t\} \cup \{B \in P \: : \: \text{\(B\) is nonouter}\}.
    \]
\end{defn}

Note that consecutive outer blocks satisfy \(\min B_{r+1}=\max B_r+1\). Merging them inserts only the unit arc \((\max B_r,\min B_{r+1})\), which cannot cross an existing arc. Hence \(\kappa(P)\) is noncrossing and has one outer block.

 \begin{exam}
Let \(P = 16|25|34|79|8\). The outer blocks of \(P\) are \(B_1=16\) and \(B_2=79\), since \(\max B_1=6<7=\min B_2\). Thus, by definition of \(\kappa\), we merge these outer blocks while leaving the nested blocks unchanged. In the arc diagram, \(\kappa\) adds the arc \textcolor{red}{\((6,7)\)}, which thereby joins the two outer blocks into a single block.
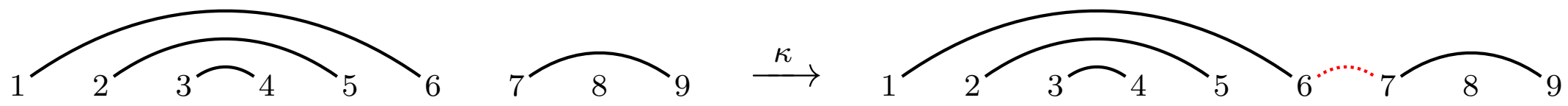
\[
\begin{tikzpicture}[scale = 0.75, baseline = 0.75*-0.2]
    \foreach \x in {1, ..., 9}
        {\draw[fill] (\x - 1, 0) node[inner sep=1pt] (\x) {$\scriptstyle \x$};}

    \foreach \i\j in {1/6, 2/5, 3/4, 7/9}
        {\draw[thick] (\i) to[out = 35, in = 145] (\j);}
\end{tikzpicture}
\quad\xrightarrow{\ \kappa\ }\quad
\begin{tikzpicture}[scale = 0.75, baseline = 0.75*-0.2]
    \foreach \x in {1, ..., 9}
        {\draw[fill] (\x - 1, 0) node[inner sep=1pt] (\x) {$\scriptstyle \x$};}

    \foreach \i\j in {1/6, 2/5, 3/4, 7/9}
        {\draw[thick] (\i) to[out = 35, in = 145] (\j);}

    \draw[thick, densely dotted, red] (6) to[out = 35, in = 145] (7);
\end{tikzpicture}
\]
\end{exam}

\begin{defn}
    Let \(P \in \operatorname{NC}_n\), and let \(B \in P\) be nonouter. Consider the blocks \(C\in P\) such that \(B\) is nested inside \(C\). This collection is nonempty and finite. Moreover, any two such blocks are nested: since both of their intervals contain the interval of \(B\), the two intervals overlap, and the noncrossing condition implies that one interval must be contained in the other. Hence, these blocks are totally ordered by nesting. Since the collection is finite, it has a unique minimal element, and we call it the \emph{immediate parent} of \(B\). Equivalently, it is the unique block \(A\in P\) such that \(B\) is nested inside \(A\) and there is no block \(C\in P\) with \(B\) nested inside \(C\) and \(C\) nested inside \(A\).

    Noncrossingness forces the entire interval of \(B\) to lie in a single gap between two consecutive elements of \(A\). Otherwise an arc of \(B\) would cross an arc of \(A\). Thus, there are unique consecutive elements \(a<a'\) of \(A\) such that
    \(
        a<\min B\le \max B<a'.
    \)
    Define
    \(
        s(B):=a'=\min\{c\in A:c>\max B\}
    \)
    and define the \emph{connector edge} of \(B\) by
    \[
        e_{\mathrm{con}}(B):=(\max B,s(B)).
    \]
    Thus every nonouter block has exactly one connector to its immediate parent.
\end{defn}

{
\begin{exam}\label{ex:connector}
Consider \(P=189|25|34|67\). The arc \((1,8)\) of the outer block
\(189\) spans two immediate child blocks, \(25\) and \(67\).
Their intervals are disjoint. One level deeper, the block \(34\)
lies under the arc \((2,5)\).
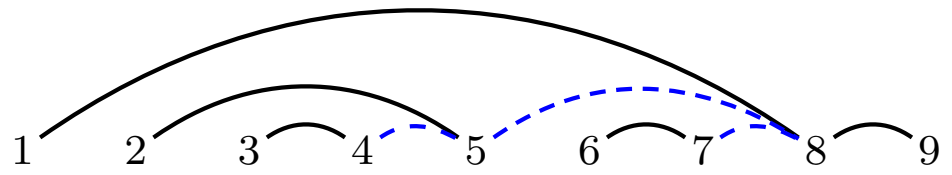
\[
\begin{tikzpicture}[scale=0.75,color=black]
    \foreach \i in {1,...,9} {
        \node[inner sep=1pt] (v\i) at (\i-1,0) {$\scriptstyle\i$};
    }
    \foreach \i/\j in {1/8,8/9,2/5,3/4,6/7} {
        \draw[thick] (v\i) to[out=35,in=145] (v\j);
    }
    \draw[blue,thick,densely dashed] (v5) to[out=35,in=148] (v8);
    \foreach \i/\j in {4/5,7/8} {
        \draw[blue,thick,dashed] (v\i) to[out=35,in=150] (v\j);
    }
\end{tikzpicture}
\]
The dashed blue arcs above the vertices show the connectors.
\end{exam}
}

\begin{defn}[Graph of a partition]\label{def:partition-graph}
The graph \(G(P)\) of a noncrossing partition \(P\) has vertex set \([n]\) and edge set
\[
    E(P)=\bigcup_{B\in P} E_{\mathrm{arc}}(B)
       \ \cup\ \{e_{\mathrm{con}}(B):B\in P\text{ is nonouter}\}.
\]
Thus \(G(P)\) contains the consecutive arcs inside every block and one connector joining every nonouter block to its immediate parent.
\end{defn}

\begin{lem}\label{lem:graph-forest}
For every \(P\in\operatorname{NC}_n\), the graph \(G(P)\) is a graph-theoretic forest with exactly one connected component for each outer block of \(P\). In particular, if \(o(P)\) is the number of outer blocks, then
\[
    |E(P)|=n-o(P).
\]
\end{lem}
\begin{proof}
The arcs of each block form a path. Contracting each of these paths to a vertex leaves just the connectors, one from each nonouter block to its immediate parent. Following parent edges leads to a unique outer block. Thus the components are indexed by the outer blocks. The block paths contribute \(n-|P|\) edges and the connectors contribute \(|P|-o(P)\), giving \(n-o(P)\) edges in total. A graph on \(n\) vertices with \(o(P)\) components and \(n-o(P)\) edges is acyclic.
\end{proof}

{
\begin{defn}[The cube of a partition]\label{def:cube-coordinates}
For \(P\in\operatorname{NC}_n\), let
\[
    Q_P=[0,1]^{E(P)}.
\]
At a vertex of \(Q_P\), an edge \(e\) is \emph{selected} if \(x_e=0\)
and \emph{omitted} if \(x_e=1\). For a set \(S\subseteq E(P)\) of
selected edges, let \(\pi_P(S)\) be the partition of \([n]\) into the
connected components of \(([n],S)\). This partition labels the vertex.
\end{defn}
}

These vertex labels are noncrossing. To see this, draw all the edges
of \(G(P)\) above the line. Within each gap of a parent block, the
intervals of its immediate children are disjoint. A connector runs from
the last element of a child to the right endpoint of the gap, passing
over the entire intervals of any intervening children. Connectors in the
same gap share their right endpoint. Thus the edges can be drawn without
crossings. If two components of a subgraph had interleaved vertices, paths
in the two components joining the alternating endpoints would have to
intersect. Moreover, the label determines \(S\) within \(Q_P\). By
Lemma~\ref{lem:graph-forest}, an edge of \(G(P)\) is selected exactly
when its endpoints belong to the same block of \(\pi_P(S)\).

{
\begin{exam*}
In Example~\ref{ex:connector}, select all the block arcs and omit all
three connectors. The vertex label is \(P=189|25|34|67\). Selecting the
connector \((4,5)\) joins the blocks \(34\) and \(25\), giving the label
\(189|2345|67\). These two labels are the endpoints of the edge of \(Q_P\)
on which only \(x_{(4,5)}\) varies.
\end{exam*}
}

To assemble the cubes, associate to every nonempty coordinate face its set of
vertex labels. Explicitly, let \(A\) consist of the edges fixed to be
selected, and let \(B\) consist of these edges together with all edges
whose coordinates vary. The vertex labels on this face are precisely
\[
    \{\pi_P(S): A\subseteq S\subseteq B\}.
\]
Identify faces from different cubes when their entire sets of
vertex labels agree. Order the resulting sets by inclusion and adjoin the
empty face. The following proposition justifies this description of
\(\mathcal C_n\).

{
\begin{prop}[Partition model of the moment complex]\label{prop:partition-model}
The poset just defined is isomorphic, by matching vertex labels, to the
face poset of \(\operatorname{Complex}(\operatorname{QFl}_n)\), including the empty face. The
Bruhat maximum of the vertex labels of \(Q_P\) is \(P\). More generally,
every nonempty face \(H\) has a unique Bruhat-maximum vertex label
\(\lambda(H)\), and \(H\) is a face of \(Q_{\lambda(H)}\).
\end{prop}
}
We prove Proposition~\ref{prop:partition-model} in
Appendix~\ref{app:comparison}. In what follows, a face of \(Q_P\) and
its image in \(\mathcal C_n\) are denoted by the same symbol. All cubes and faces used in unions are closed cells.

\section{Proof of the First-Carrier Theorem}
 
For \(P\in\operatorname{NC}_n\) and \(j\in[n]\), define
\[
    d_P(j):=\#\{B\in P:\min B<j<\max B,\ j\notin B\}.
\]
Thus, \(d_P(j)\) is the number of blocks of \(P\) that strictly nest \(j\). Define the vertical code
\[
    \operatorname{vert}_j(P)=
    \begin{cases}
        2d_P(j),&j\text{ is maximal in its block},\\
        2d_P(j)+1,&j\text{ is not maximal in its block}.
    \end{cases}
\]
Gobet--Williams \cite[Theorem~1.2 and Section~6.1]{GobetWilliams2016} give the following criterion for Bruhat order. We use its formulation in Proposition~8.6 of \cite{BergeronGagnonNadeauSpinkTewari2025}, with the decreasing-cycle convention specified above. For \(P,Q \in \operatorname{NC}_n\),
\[
    P\le_B Q
    \quad\Longleftrightarrow\quad
    \operatorname{vert}_j(P)\le \operatorname{vert}_j(Q)\text{ for all }j.
\]

\begin{lem}\label{lem:kappa}
    For all \(P\in\operatorname{NC}_n\):
    \begin{enumerate}
        \item[(i)] \(P \leq_B \kappa(P)\),
        \item[(ii)] \(\kappa \) is Bruhat-order preserving,
        \item[(iii)] \(\kappa(P) \) is the least fully supported Bruhat upper bound of \(P\).
    \end{enumerate}
\end{lem}
\begin{proof}
    Let \(B_1, \dots, B_t\) be the outer blocks of \(P\), indexed as in Definition~\ref{def:kappa}. Applying \(\kappa\) does not change \(d_P(j)\). If \(j\) lies in an outer block, neither that block nor the merged outer block nests \(j\). On the other hand, if \(j\) lies in a nonouter block, its unique outer ancestor is replaced by the merged outer block, while all other nesting relations are unchanged. Thus merging the outer blocks changes the block-maximal status only at 
    \(
        \max B_1, \dots,  \max B_{t-1}.
    \)
    Each changes from maximal to non-maximal. Since \(\max B_t=n\), it follows that
\[
    \operatorname{vert}_j(\kappa(P))=
    \begin{cases}
        \max\{1,\operatorname{vert}_j(P)\},&j<n,\\
        0,&j=n.
    \end{cases}
\]
Hence \(\operatorname{vert}_j(P)\le \operatorname{vert}_j(\kappa(P))\) for every \(j\), so Proposition~8.6 of \cite{BergeronGagnonNadeauSpinkTewari2025} gives \(P\le_B\kappa(P)\).

If \(P\leq_B Q\), then \(\operatorname{vert}_j(P)\le\operatorname{vert}_j(Q)\) for every \(j\). Applying the monotone map \(a\mapsto\max\{1,a\}\) for \(j<n\) gives
\(
    \operatorname{vert}_j(\kappa(P))\le\operatorname{vert}_j(\kappa(Q))
\)
for every \(j\), and therefore \(\kappa(P)\le_B\kappa(Q)\). This proves (ii).

If \(P\le_B\nu\) and \(\nu\) is fully supported, then \(\kappa(\nu)=\nu\), so (ii) gives \(\kappa(P)\le_B\kappa(\nu)=\nu\). Since \(\kappa(P)\) is fully supported and lies above \(P\) by (i), this proves (iii).
\end{proof}
 

\begin{lem}[Boundary-coordinate face lemma]\label{lem:boundary-face}
    Let \(P\in\operatorname{NC}_n\), and let \(B_1,\ldots,B_t\) be its outer blocks, indexed as in Definition~\ref{def:kappa}. For \(1\le r<t\), set
    \(
        e_r:=(\max B_r,\min B_{r+1}).
    \)
    Then \(e_r\in E(\kappa(P))\), and
    \[
        Q_P
        =\bigcap_{r=1}^{t-1}
            \{x\in Q_{\kappa(P)}:x_{e_r}=1\}.
    \]
    When \(t=1\), the intersection over \(1\le r<t\) is understood to be all of \(Q_{\kappa(P)}\).
\end{lem}
\begin{proof}
    The intervals spanned by the outer blocks cover \([n]\) in the specified order, so
    \(
        \min B_{r+1}=\max B_r+1.
    \)
    Merging these blocks adds precisely the arcs \(e_1,\ldots,e_{t-1}\).
    Every connector keeps its endpoints. Indeed, the only parent blocks changed
    by the merger are outer, and each child of such a block remains in the
    same gap between two consecutive elements of that block. Hence
    \[
        E(\kappa(P))=E(P)\sqcup\{e_1,\ldots,e_{t-1}\}.
    \]
    On the indicated face of \(Q_{\kappa(P)}\) all the added edges
    are omitted. The remaining selections are exactly the subsets
    \(S\subseteq E(P)\). Each gives the same graph \(([n],S)\), and
    therefore the same vertex label, in both cubes. The face is identified
    with \(Q_P\) by the gluing rule.
\end{proof}

\begin{exam}
Let \(P=12|34\).
Both blocks are outer, so \(\kappa(P)=1234\).
The graph \(G(P)\) has edges \((1,2)\) and \((3,4)\).
Merging the blocks adds the edge \((2,3)\). Hence
\(
    Q_P=\{x\in Q_{\kappa(P)}:\ x_{(2,3)}=1\}\subseteq Q_{\kappa(P)},
\)
so \(Q_P\) is a two-dimensional facet of the three-cube \(Q_{\kappa(P)}\).
\end{exam}

\begin{rem}[Top-dimensional labels]\label{rem:top-labels}
By Lemma~\ref{lem:graph-forest}, \(\dim Q_P=n-o(P)\). Thus
\(Q_P\) has dimension \(n-1\) exactly when \(P\) is fully supported.
By Lemma~\ref{lem:boundary-face} every other \(Q_P\) is a proper face
of \(Q_{\kappa(P)}\). Consequently the maximal cells of \(\mathcal C_n\)
are precisely the cubes indexed by fully supported partitions.
\end{rem}

\begin{thm}\label{thm:first-carrier}
    Let \(H\) be any nonempty global face, and let \(P = \lambda(H)\). The first top cube carrying \(H\), in any Bruhat-compatible top-cube order, is \(Q_{\kappa(P)}\).
\end{thm}
\begin{proof}
    Proposition~\ref{prop:partition-model} and Lemma~\ref{lem:boundary-face}
    give \(H\subseteq Q_P\subseteq Q_{\kappa(P)}\). Now let \(Q_\nu\) be
    any maximal cube containing \(H\). Its vertex labels include \(P\)
    and have Bruhat maximum \(\nu\), so \(P\le_B\nu\). By
    Remark~\ref{rem:top-labels}, \(\nu\) is fully supported. Lemma~\ref{lem:kappa}(iii)
    therefore gives \(\kappa(P)\le_B\nu\). Thus \(Q_{\kappa(P)}\) has
    the least label among the maximal cubes carrying \(H\), and comes
    first among them in every linear extension.
\end{proof}

\section{The new vertices of a maximal cube}\label{sec:attaching}
Fix a fully supported partition \(\mu\in\operatorname{NC}_n\), and let \(B_0\) be its unique outer block. The edges of \(G(\mu)\) index the coordinates of \(Q_\mu\).

\begin{defn}[Essential, connector, and free coordinates]\label{def:edge-types}
Define
\[
\begin{aligned}
E_{\mathrm{ess}}(\mu)
   := \bigcup_{\substack{B\in\mu\\B\neq B_0}}E_{\mathrm{arc}}(B) \cup
   \{(a,b)\in E_{\mathrm{arc}}(B_0):b>a+1\},
\end{aligned}
\]
\[
    E_{\mathrm{con}}(\mu)
      :=\{e_{\mathrm{con}}(B):B\in\mu,\ B\neq B_0\},
\]
and
\[
    E_{\mathrm{free}}(\mu)
      :=E(\mu)\setminus
        \bigl(E_{\mathrm{ess}}(\mu)\cup E_{\mathrm{con}}(\mu)\bigr).
\]
Thus, \(E_{\mathrm{free}}(\mu)\) consists exactly of the unit arcs \((a,a+1)\) between consecutive elements of the outer block \(B_0\).
\end{defn}

\begin{exam}\label{ex:edge-types}
Consider
\(
    \mu=14679|23|5|8.
\)
The essential edges are
\(
    \textcolor{darkgreen}{(1,4)},\textcolor{darkgreen}{(4,6)},\textcolor{darkgreen}{(7,9)},\textcolor{darkgreen}{(2,3)};
\)
the connector edges are
\(
    \textcolor{blue}{(3,4)},\textcolor{blue}{(5,6)},\textcolor{blue}{(8,9)};
\)
and the free edge is \textcolor{red}{\((6,7)\)}.

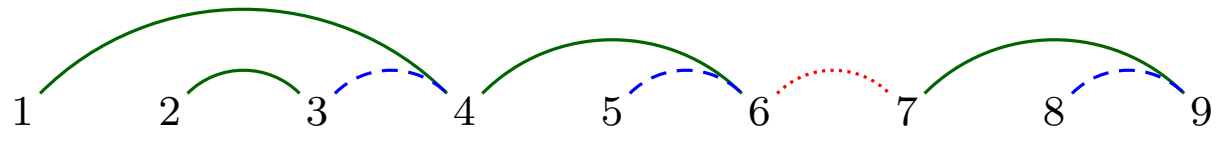
\begin{center}
\begin{tikzpicture}[scale=1.5, baseline=-0.5ex]
    \foreach \x in {1,2,...,9}{
        \node[inner sep=1pt] (\x) at ({0.65*(\x-1)},0) {\scriptsize \x};
    }

    \draw[darkgreen, semithick] (1) to[out=45,in=135] (4);
    \draw[darkgreen, semithick] (4) to[out=45,in=135] (6);
    \draw[darkgreen, semithick] (2) to[out=45,in=135] (3);
    \draw[darkgreen, semithick] (7) to[out=45,in=135] (9);

    \draw[red, densely dotted, semithick] (6) to[out=45,in=135] (7);

    \draw[blue, densely dashed, semithick] (3) to[out=45,in=135] (4);
    \draw[blue, densely dashed, semithick] (5) to[out=45,in=135] (6);
    \draw[blue, densely dashed, semithick] (8) to[out=45,in=135] (9);
\end{tikzpicture}
\end{center}
Green denotes essential edges, dashed blue denotes connector edges, and dotted red  denotes free edges.
\end{exam}

For a vertex \(v\in Q_\mu\), write \(S(v)=\{e\in E(\mu):x_e(v)=0\}\).
Its label is \(\pi_\mu(S(v))\) by Definition~\ref{def:cube-coordinates}.

\begin{lem}\label{lem:new}
    For a vertex with selected-edge set \(S\subseteq E(\mu)\), 
    \[
        \kappa(\pi_\mu(S)) = \mu \qquad \Longleftrightarrow \qquad \text{\(E_{\text{ess}}(\mu) \subseteq S\) and \(E_{\text{con}}(\mu) \cap S = \varnothing\)}.
    \]
\end{lem}
\begin{proof}
    Put \(P=\pi_\mu(S)\). First, suppose that \(\kappa(P)=\mu\).

    If a connector edge \(e_{\mathrm{con}}(B)\) were selected, then one connected component of \(([n],S)\) would contain elements from two distinct blocks of \(\mu\), namely \(B\) and its immediate parent. Since \(\kappa\) only merges blocks and never separates them, those elements would remain in one block of \(\kappa(P)\), and \(\kappa(P) \neq \mu\). Thus, no connector can be selected, and \(E_{\mathrm{con}}(\mu)\cap S=\varnothing\).

    Now let \((a,b)\) be an arc of a nonouter block \(B\). If \((a,b)\notin S\), then the path of consecutive arc edges through \(B\) is broken. Every edge of \(G(\mu)\) with exactly one endpoint in \(B\) is a connector edge: either the connector of \(B\) to its parent or the connector of one of its children to \(B\). As we have shown, no connector edge is selected. Thus, no selected path can leave \(B\) and later re-enter it. Hence, the vertices of \(B\) lie in at least two distinct connected components of \(([n],S)\), so \(B\) is split between at least two blocks of \(P\). Applying \(\kappa\) cannot reconstruct \(B\). Indeed, nonouter pieces are never merged with one another, while every piece that is outer is merged with all outer blocks of \(P\), including an outer component containing vertices of the outer block \(B_0\) of \(\mu\). Thus, an outer piece of \(B\) is merged with vertices outside \(B\). In all cases \(B\) fails to occur as a block of \(\kappa(P)\), contradicting \(\kappa(P)=\mu\). Hence, every arc of every nonouter block belongs to \(S\).

    Finally, let \((a,b)\) be an outer arc of \(B_0\) with \(b>a+1\). Since \(a\) and \(b\) are consecutive elements of \(B_0\), every element of the open interval \((a,b)\) belongs to a block contained in that gap. Otherwise an arc would cross \((a,b)\). In particular, the gap contains at least one block. Choose a block \(C\) contained in the gap that is not nested inside any other block contained in the gap. Then no block lies strictly between \(C\) and \(B_0\) in the nesting order, so \(C\) is an immediate child of \(B_0\). If \((a,b)\notin S\), then the selected outer-arc path of \(B_0\) is cut at this gap. Since \(a\) and \(b\) are consecutive in \(B_0\), the component of \(B_0\) containing \(a\) has maximum \(a<\min C\), while the component containing \(b\) has minimum \(b>\max C\). Thus no component coming from \(B_0\) can contain the interval of \(C\). By the choice of \(C\), no block lying in the gap properly contains \(C\), and because no connector edge is selected, distinct blocks of \(\mu\) cannot become joined in \(([n],S)\). Hence no block of \(P\) has interval strictly containing the interval of \(C\), so \(C\) is an outer block of \(P\). Applying \(\kappa\) then merges \(C\) with the other outer blocks, whereas \(C\) is a distinct nonouter block of \(\mu\). This contradicts \(\kappa(P)=\mu\). Therefore every long outer arc belongs to \(S\). This proves necessity.

    Conversely, suppose \(E_{\text{ess}}(\mu) \subseteq S\) and  \(E_{\text{con}}(\mu) \cap S = \varnothing\). Every nonouter block is a block of \(P\), since all its arcs are selected and no connector is selected. The only outer arcs that may be omitted are free edges \((a,a+1)\), so they can only split \(B_0\) into adjacent components, with one ending at \(a\) and the next beginning at \(a+1\) for each omitted arc \((a,a+1)\). We claim that every nonouter block remains nonouter in \(P\). Let \(B\) be nonouter and let \(A\) be its immediate parent. If \(A\neq B_0\), then all consecutive arc edges of \(A\) are selected, so the block of \(P\) containing \(A\) still spans \(B\). If \(A=B_0\), then \(B\) lies in a gap between consecutive elements \(a<b\) of \(B_0\). Such a gap contains \(B\), so \(b>a+1\); hence the outer arc \((a,b)\) is essential and is selected. Thus the component of \(B_0\) containing both \(a\) and \(b\) still spans \(B\). In either case, \(B\) is nested inside a block of \(P\). Iterating this argument up the finite parent chain shows that every nonouter block of \(\mu\) is nested inside a component of \(B_0\) in \(P\). Therefore the outer blocks of \(P\) are exactly the components into which \(B_0\) was split. Applying \(\kappa\) merges precisely these components and leaves every nonouter block unchanged, yielding \(\mu\). Hence \(\kappa(P)=\mu\).
\end{proof} 

\begin{rem}\label{rem:new-vertex-face}
Define
\[
    S_\mu:=E_{\mathrm{ess}}(\mu)\sqcup E_{\mathrm{con}}(\mu),  \qquad
    \varepsilon_e:=
    \begin{cases}
    0,&e\in E_{\mathrm{ess}}(\mu),\\
    1,&e\in E_{\mathrm{con}}(\mu).
    \end{cases}
\]
Let
\[
R_\mu:=\{x\in Q_\mu:\ x_e=\varepsilon_e\ \text{for every }e\in S_\mu\}\subseteq Q_\mu.
\]
Thus \(R_\mu\) is the coordinate face obtained by fixing the
essential-edge coordinates to \(0\) and the connector coordinates to \(1\),
while the coordinates in \(E_{\mathrm{free}}(\mu)\) vary independently
in \([0,1]\). In Example~\ref{ex:edge-types}, the only varying coordinate
is \(x_{(6,7)}\), so \(R_\mu\) is an edge.
For a vertex \(v\in Q_\mu\), Lemma~\ref{lem:new} gives
\[
v\in R_\mu \quad\Longleftrightarrow\quad \kappa(\pi_\mu(S(v)))=\mu.
\]
By Theorem~\ref{thm:first-carrier-intro}, this is equivalent to saying that \(Q_\mu\) is the first top cube carrying \(v\). Thus, the vertices appearing for the first time in \(Q_\mu\) are exactly the vertices of \(R_\mu\).
\end{rem}

\section{Attaching loci and shellability}

For a cube \(Q\) with coordinates \((x_e)_{e\in I}\), we write
\[
    \{x_e=a\}_Q:=\{x\in Q:x_e=a\}
    \qquad(e\in I,\ a\in\{0,1\}).
\]
More generally, for \(J\subseteq I\) and \(a_e\in\{0,1\}\), we set
\[
    \{x_e=a_e:e\in J\}_Q
    :=\{x\in Q:x_e=a_e\text{ for every }e\in J\}.
\]
The subscript specifies the cube in whose coordinates the equations
are interpreted. These sets are coordinate faces of \(Q\), with all
unspecified coordinates allowed to vary in \([0,1]\).

\begin{thm}\label{thm:attaching}
Let \(Q_1,\ldots,Q_N\) be the top cubes in any Bruhat-compatible
linear extension, and let \(\mu=\mu_{j}\) be the fully supported label
of \(Q_{j}\). With \(S_\mu\) and \(\varepsilon_e\) as in
Remark~\ref{rem:new-vertex-face}, we have
\[
    Q_{j}\cap Y_{j-1} =\bigcup_{e\in S_\mu}\{x_e=1-\varepsilon_e\}_{Q_j}.
\]
\end{thm}
\begin{proof}
Fix \(j\), put \(\mu=\mu_{j}\), and write
\[
    L_\mu :=\bigcup_{e\in S_\mu}\{x_e=1-\varepsilon_e\}_{Q_j}.
\]
We prove that \(Q_{j}\cap Y_{j-1}=L_\mu\).

Let \(G\) be a nonempty coordinate face of \(Q_{j}\) contained in \(L_\mu\).
Let \(P=\lambda(G)\) and let \(v_P\in G\) be the vertex with label
\(P\). Since \(v_P\in G\subseteq L_\mu\), there exists \(e\in S_\mu\) such that
\(
    x_e(v_P)=1-\varepsilon_e.
\)
Thus \(v_P\notin R_\mu\), and Remark~\ref{rem:new-vertex-face} gives
\(\kappa(P)\ne\mu\).
Since \(P\) is a vertex label of \(Q_{j}\), we have \(P\le_B\mu\).
The label \(\mu\) is fully supported, so Lemma~\ref{lem:kappa}(iii) gives \(\kappa(P)<_B\mu\), with strict inequality. By the first-carrier theorem, the first top cube containing \(G\) is \(Q_{\kappa(P)}\). Its label occurs earlier in the chosen linear extension, so there exists \(i<j\) such that
\[
    G\subseteq Q_{\kappa(P)}=Q_{i}\subseteq Y_{j-1}.
\]
As \(L_\mu\) is a union of coordinate facets of \(Q_{j}\), this
proves \(L_\mu\subseteq Q_{j}\cap Y_{j-1}\).

Conversely, let \(H\) be a coordinate face of \(Q_{j}\) not contained
in \(L_\mu\). For each \(e\in S_\mu\), the coordinate \(x_e\)
is either fixed to \(\varepsilon_e\) on \(H\) or is allowed to vary: it cannot
be fixed to \(1-\varepsilon_e\). Set all varying \(S_\mu\)-coordinates
to \(\varepsilon_e\), and assign either \(0\) or \(1\) to each of the remaining unrestricted coordinates. This gives a vertex \(v\in H\cap R_\mu\). By Remark~\ref{rem:new-vertex-face}, the first top cube containing
\(v\) is \(Q_{j}\). Therefore \(v\notin Y_{j-1}\), so
\(H\not\subseteq Y_{j-1}\).

The intersection \(Q_{j}\cap Y_{j-1}\) is a subcomplex of the cubical
cell \(Q_{j}\), hence a union of coordinate faces. By the preceding
argument every one of these faces is contained in \(L_\mu\).
Thus \(Q_{j}\cap Y_{j-1}\subseteq L_\mu\), proving equality.
\end{proof}

\begin{exam} Consider the attachment for \(n=3\), with the vertex labels shown in Figure~\ref{fig:n3-common-path}. There are two fully supported noncrossing partitions of \([3]\), namely \(\mu_1=123\) and \(\mu_2=13|2\), in this particular Bruhat order. The top cells indexed by \(\mu_1\) and \(\mu_2\) are \(Q_1\) and \(Q_2\), respectively, both of which are squares. For \(\mu_1\) the essential and connector sets are empty, so \(Q_1\) is attached first with empty intersection. For \(\mu_2\), the unique long outer arc \((1,3)\) is essential and \((2,3)\) is a connector, so
\(S_{\mu_2}=\{(1,3),(2,3)\}\), with
\(\varepsilon_{(1,3)}=0\) and \(\varepsilon_{(2,3)}=1\). Hence
\[
    Q_2\cap Q_1=\bigcup_{e\in S_{\mu_2}}\{x_e=1-\varepsilon_e\}_{Q_2}.
\]
The two facets in this union are adjacent edges of the square \(Q_2\). Their union is a closed star of a vertex in \(\partial Q_2\), but is not a single face of that square. In particular, \(\mathcal C_3\) is not a polytopal complex in the usual sense, that is, the intersection of two cells is not necessarily a common face of both.
\end{exam}

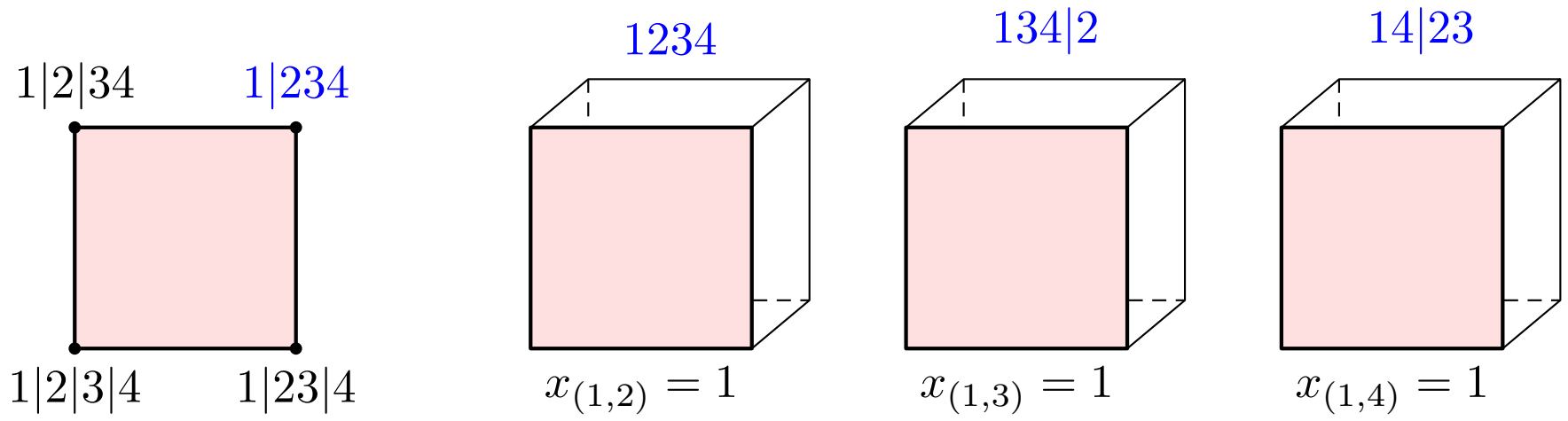
\begin{figure}[htbp]
\centering
\begin{tikzpicture}[font=\small, draw=black, line join=round]
  \begin{scope}
    \filldraw[fill=red!12,thick] (0,0) rectangle (1.65,1.65);
    \foreach \x/\y in {0/0,1.65/0,0/1.65,1.65/1.65}
      \fill (\x,\y) circle (1.3pt);
    \node[below] at (0,0) {$1|2|3|4$};
    \node[below] at (1.65,0) {$1|23|4$};
    \node[above] at (0,1.65) {$1|2|34$};
    \node[above] at (1.65,1.65) {\textcolor{blue}{$1|234$}};
  \end{scope}
  \foreach \offset/\labeltext/\j in
      {3.4/{1234}/2, 6.2/{134|2}/3, 9.0/{14|23}/4}{
    \begin{scope}[xshift=\offset cm]
      \draw[densely dashed]
        (0,0) -- (.43,.36) -- (.43,2.01)
        (.43,.36) -- (2.08,.36);
      \draw (0,1.65) -- (.43,2.01) -- (2.08,2.01)
        -- (2.08,.36) -- (1.65,0)
        (1.65,1.65) -- (2.08,2.01);
      \filldraw[fill=red!12,thick] (0,0) rectangle (1.65,1.65);
      \node[above] at (1.04,2.06) {\textcolor{blue}{$\labeltext$}};
      \node[below] at (.825,0) {$x_{(1,\j)}=1$};
    \end{scope}
  }
\end{tikzpicture}
\caption{A common square in \(\mathcal C_4\). The three shaded facets
are identified with the labeled square on the left. In all four drawings,
moving right selects \((2,3)\) and moving up selects \((3,4)\) on the
shaded face. The label above each cube is its fully supported partition.}
\label{fig:n4-common-facet}
\end{figure}

{
\begin{exam}\label{ex:three-carriers}
For \(n=4\), a single square can be a facet of three maximal cubes.
Consider the square in Figure~\ref{fig:n4-common-facet}. Its two varying
coordinates correspond to \((2,3)\) and \((3,4)\). Any cube \(Q_\mu\) containing
the square must have both edges in its graph \(G(\mu)\). Indeed, the
labels \(1|23|4\) and \(1|2|34\) each contain a two-element block,
whose two vertices must be joined by a selected edge.
Of the five top labels, exactly the three shown have both edges. In
each of these cubes the third edge is omitted throughout the square.
Thus the three cubes meet along the interior of a common facet. Every triangle in the barycentric subdivision of this square is incident to three tetrahedra in \(\operatorname{sd}(\mathcal C_4)\). The link of such a triangle therefore consists of three points, so \(\mathcal C_4\) is not a manifold with boundary.
\end{exam}
}

For a nonempty face \(F\) of the boundary of a cube \(Q\),
we use the closed-star convention
\[
    \operatorname{st}_{\partial Q}(F)
    :=\bigcup_{\substack{D\text{ a facet of }Q\\F\subseteq D}}D.
\]
Thus its facets contain \(F\). In the following lemma, the attaching
facets are disjoint from the new-vertex face \(R_\mu\) and contain its
opposite face \(R_\mu^-\).

A \(d\)-dimensional \emph{piecewise linear (PL) ball} is a polyhedron
admitting a homeomorphism to a closed \(d\)-simplex that is affine on each
simplex of a suitable triangulation. For \(n=3\), the two adjacent edges
forming \(Q_2\cap Q_1\) in Figure~\ref{fig:n3-common-path} can be
straightened to a closed interval by such a map. Their union is therefore
a PL \(1\)-ball.

We write \(B^r\) for the standard closed PL \(r\)-ball \([0,1]^r\).

\begin{lem}\label{lem:ball}
Let \(d=n-1\). Set \(A_j:=Q_j\cap Y_{j-1}\) for \(1\le j\le N\). For \(j>1\), the attaching locus
is a pure \((d-1)\)-dimensional PL ball. More precisely, if
    \[
    R^-_{\mu_{j}}:=\{x_e=1-\varepsilon_e:e\in S_{\mu_{j}}\}_{Q_j},
    \]
then
    \[
    A_{j}=\operatorname{st}_{\partial Q_{j}}(R^-_{\mu_{j}}) \cong_{\mathrm{PL}}B^{d-1}.
    \]
Further, \(A_1=\varnothing\).
\end{lem}

\begin{proof}
    The equality \(A_1=\varnothing\) follows immediately from \(Y_0= \varnothing\). Now fix \( j>1 \), set \(\mu=\mu_{j}\), and let \(s=|S_\mu|\). By Lemma~\ref{lem:kappa} applied to \(1|2|\cdots|n\),
    the least fully supported label is \(\mu_1=1\cdots n\).
    Since \(j>1\), \(\mu\ne\mu_1\), so its unique outer block has an
    arc \((a,b)\) with \(b>a+1\). A block in the gap \((a,b)\) is
    nonouter and contributes a connector to \(S_\mu\). Hence \(s\ge1\).

  By Theorem~\ref{thm:attaching},
\[
A_{j}=\bigcup_{e\in S_\mu}\{x_e=1-\varepsilon_e\}_{Q_j}.
\]
By Remark~\ref{rem:new-vertex-face}, the new vertices are the vertices of \(R_\mu\). The facets in this union are precisely those disjoint from \(R_\mu\), proving Theorem~\ref{thm:attaching-intro}.
The facets of \(Q_{j}\) containing the face
\[
R^-_\mu=\{x_e=1-\varepsilon_e:e\in S_\mu\}_{Q_j}
\]
are exactly the facets appearing in this union. Hence, by definition of the closed star,
\(
A_{j}=\operatorname{st}_{\partial Q_{j}}(R^-_\mu).
\)

It remains to show that this closed star is PL homeomorphic to a \((d-1)\)-ball. After permuting coordinates and reflecting some coordinates, we may identify \(Q_{j}\) with \([0,1]^d\) so that 
\[R^-_\mu=\{x_1=\cdots=x_s=0\}_{[0,1]^d} .\]
Then
    \[
    A_{j}=  \left( \bigcup_{k=1}^s \{x_{k}=0\}_{[0,1]^s} \right) \times[0,1]^{d-s}.
    \]
Now, set
        \[
        K_s:= \bigcup_{k=1}^s \{x_{k}=0\}_{[0,1]^s}.
        \]
We claim that \( K_s\cong_{\mathrm{PL}}B^{s-1} \).
We prove this by induction on \(s\). For \(s=1\), \(K_1\) is a point, hence a PL \(0\)-ball. For \(s>1\), write a point of \([0,1]^s\) as \((x_1, \dots, x_{s-1},x_s)\). We know the set \(K_s \) consists of points for which at least one coordinate is \(0\). There are now two possibilities. First, one of the first \(s-1\) coordinates is \(0\). Those points form \(K_{s-1}\times[0,1]\). Second, the last coordinate is \(0\). Those points form \( [0,1]^{s-1}\times\{0\}\). Thus, 
        \[
        K_s = (K_{s-1}\times[0,1]) \cup ([0,1]^{s-1}\times\{0\}).
        \]
By induction,
        \[
        K_{s-1}\times[0,1] \cong_{\mathrm{PL}} B^{s-2}\times B^1 \cong_{\mathrm{PL}} B^{s-1},
        \]
while
    \[
    [0,1]^{s-1}\times\{0\} \cong_{\mathrm{PL}}B^{s-1}.
    \]
Their intersection is \(K_{s-1}\times\{0\}\) and by the inductive assumption, 
    \[
    K_{s-1}\times\{0\} \cong_{\mathrm{PL}}B^{s-2},
    \]
and this is a PL ball contained in the boundary of each of the two \((s-1)\)-balls. Gluing two PL \((s-1)\)-balls along a common PL \((s-2)\)-ball contained in their boundaries yields a PL \((s-1)\)-ball~\cite[Section~12.2(2)]{Bjorner1995}. Thus
    \(
    K_s\cong_{\mathrm{PL}}B^{s-1},
    \)
and therefore
    \[ 
    A_{j} \cong_{\mathrm{PL}} K_s\times[0,1]^{d-s} \cong_{\mathrm{PL}} B^{s-1}\times B^{d-s} \cong_{\mathrm{PL}} B^{d-1}.
    \]
\end{proof}

We use Bj\"orner's recursive definition of shelling for a pure regular CW complex \cite[Definition~4.1]{Bjorner1984}. In dimension zero any ordering of the maximal cells is a shelling. In positive dimension the boundary of the first maximal cell must be shellable. Each subsequent cell must meet the preceding union in a pure subcomplex of codimension one whose facets come first in some shelling of the boundary of that cell.

\begin{lem}[Closed stars in cube boundaries]\label{lem:cube-star-shelling}
Let \(Q\) be a \(d\)-dimensional cube with \(d\ge1\), and let \(R\) be a nonempty face of \(\partial Q\). Then \(\operatorname{st}_{\partial Q}(R)\) is an initial segment of a shelling of \(\partial Q\).
\end{lem}

\begin{proof}
After identifying \(Q\) with the standard cube \([0,1]^d\), and then permuting or reflecting coordinates, write
\[
    R=\{x_1=\cdots=x_s=0\}_{Q}
\]
for some \(1\le s\le d\). Let
\[
    F_j^-:=\{x_j=0\}_{Q}, \qquad \text{and} \qquad F_j^+:=\{x_j=1\}_{Q}
\]
denote the facets of \(Q\). Now, since every point of \(R \) satisfies \(x_1=\cdots=x_s=0\) we have \(R \subseteq F_j^-\) for every \(j = 1,\dots, s\). No other facet contains \(R\). Therefore, 
\[
    \operatorname{st}_{\partial Q}(R)=F_1^-\cup\cdots\cup F_s^-.
\]
We claim that
\begin{equation}
    F_1^-,\ldots,F_d^-,F_1^+,\ldots,F_d^+ \label{eq:cube-shelling}
\end{equation}
is a shelling of \(\partial Q\). 

The proof is done by induction on \(d\). For \(d=1\), the boundary consists of two points and the assertion is immediate. Assume \(d>1\) and the inductive claim holds for \((d-1)\)-cubes.

The first facet \(F_1^-\) is a \((d-1)\)-cube, whose boundary is shellable by induction. Now fix \(2\le k\le d\). At the moment when \(F_k^-\) is added, the preceding facets are \(F_1^-,\ldots,F_{k-1}^-\). Hence
\[
F_k^-\cap \bigl(F_1^-\cup\cdots\cup F_{k-1}^-\bigr) = \bigcup_{j<k}(F_k^-\cap F_j^-).
\]

Inside \(F_k^-\), the faces \(F_k^-\cap F_j^-\) for \(j<k\) are the first \(k-1\) negative coordinate facets of a \((d-1)\)-cube. They therefore form an initial segment of the inductive shelling of \(\partial F_k^-\).

For \(1\le k\le d\), the intersection of \(F_k^+\) with the preceding facets is
\[
    \bigcup_{j\ne k}(F_k^+\cap F_j^-)
    \ \cup\ 
    \bigcup_{j<k}(F_k^+\cap F_j^+).
\]
Inside the \((d-1)\)-cube \(F_k^+\), this consists of all negative coordinate facets followed by the first positive coordinate facets, again an initial segment of the inductive shelling. Hence \eqref{eq:cube-shelling} is a shelling of \(\partial Q\). Its first \(s\) facets are exactly the closed star of \(R\), proving the lemma.
\end{proof}

\begin{proof}[Proof of Theorem~\ref{thm:shelling-intro}]
When \(n=1\), the unique ordering of the single vertex is a shelling. Assume \(n\ge2\). Let \(Q_1,\dots,Q_N\) be any linear extension of Bruhat order on the top-cell labels. The boundary of \(Q_1\) is shellable by Lemma~\ref{lem:cube-star-shelling}, or directly by the shelling \eqref{eq:cube-shelling}. For each \(j>1\), Lemma~\ref{lem:ball} and the proof of that lemma give
\[
    Q_{j}\cap Y_{j-1}
      =\operatorname{st}_{\partial Q_{j}}(R^-_{\mu_{j}})
\]
for a nonempty face \(R^-_{\mu_{j}}\subseteq \partial Q_{j}\). By Lemma~\ref{lem:cube-star-shelling}, this attaching locus is an initial segment of a shelling of \(\partial Q_{j}\). Hence, the top-cell order satisfies the regular-cell shelling condition at every stage.
\end{proof}

We use the abstract convention for simplicial complexes: \(\Delta\)
is a nonempty finite collection of finite sets, called faces, such that
\(\tau\subseteq\sigma\in\Delta\) implies \(\tau\in\Delta\).
Its vertex set is \(V=\bigcup_{\sigma\in\Delta}\sigma\), and
\(\dim\sigma=|\sigma|-1\); in particular, \(\dim\varnothing=-1\).
Unions and intersections of faces are operations on their vertex sets.

The geometric realization of \(\Delta\) is
\[
    |\Delta|:=\left\{(t_v)_{v\in V}\in\mathbb R_{\ge0}^{V}:
    \sum_{v\in V}t_v=1,\quad
    \{v\in V:t_v>0\}\in\Delta\right\}.
\]
Each nonempty face gives a filled simplex on its vertices, and these
simplices are glued along their common faces. Homology, connectivity, and
other topological properties of a simplicial complex refer to its
geometric realization; the same convention applies to links.

For a face \(\sigma\in\Delta\), its link is the abstract simplicial complex
\[
    \operatorname{lk}_\Delta(\sigma)
    :=\{\tau\in\Delta:\tau\cap\sigma=\varnothing,
                          \ \tau\cup\sigma\in\Delta\}.
\]
The condition \(\tau\cup\sigma\in\Delta\) means that the vertices
of the two disjoint faces together span a simplex of \(\Delta\).

Let \(A\) be a field or \(\mathbb Z\). Then
a pure simplicial complex \(\Delta\) is Cohen--Macaulay over \(A\) if
\[
    \widetilde H_i(\operatorname{lk}_\Delta(\sigma);A)=0
    \quad\text{for every face }\sigma
    \text{ and every }i<\dim\operatorname{lk}_\Delta(\sigma).
\]
The empty face is included, so this also imposes the corresponding
homology vanishing on \(\Delta\) itself. Over a field, this is equivalent
to the Stanley--Reisner ring being Cohen--Macaulay.

{
\begin{defn}[Homotopy Cohen--Macaulayness]
A pure \(d\)-dimensional simplicial complex \(\Delta\) is \emph{homotopy Cohen--Macaulay} if, for every face \(\sigma\), including the empty face, the link \(\operatorname{lk}_\Delta(\sigma)\) is \((d-\dim\sigma-2)\)-connected. Here \((-1)\)-connected means nonempty, and \((-2)\)-connected imposes no condition.
\end{defn}

\begin{proof}[Proof of Corollary~\ref{cor:cm-intro}]
The proof of Theorem~\ref{thm:shelling-intro} applies to each prefix \(Q_1,\ldots,Q_j\), so \(Y_j\) is a pure shellable regular CW complex of dimension \(n-1\). By \cite[Proposition~4.4(a)]{Bjorner1984}, its barycentric subdivision is pure and shellable. Pure shellability implies homotopy Cohen--Macaulayness \cite[Section~11.5]{Bjorner1995}. It also gives vanishing reduced integral homology of each link below the dimension of that link, which is the Cohen--Macaulay condition over \(\mathbb Z\). The universal coefficient theorem gives the same vanishing over every field.
\end{proof}

\begin{rem}[Double and uniform Cohen--Macaulayness]
A simplicial complex is \emph{doubly Cohen--Macaulay} if it is Cohen--Macaulay and deleting any vertex gives a Cohen--Macaulay complex of the same dimension. It is \emph{uniformly Cohen--Macaulay} if the analogous condition holds upon removing any maximal face while retaining its proper faces \cite{MatsuokaMurai}.

For every \(n\ge2\), the complex \(\Delta=\operatorname{sd}(\mathcal C_n)\)
is neither uniformly nor doubly Cohen--Macaulay over any field.
Indeed, fix a field \(\Bbbk\) and put \(d=n-1\).
The complex \(\Delta\) is contractible
\cite[Theorem~9.3 and Remark~9.17]{CoxeterFlag2026}.
If \(F\) is a maximal simplex and \(\Gamma=\Delta\setminus\{F\}\),
with all proper faces of \(F\) retained, then the relative chain complex gives
\[
    H_d(\Delta,\Gamma;\Bbbk)\cong\Bbbk.
\]
The reduced homology sequence of the pair therefore gives
\(\widetilde H_{d-1}(\Gamma;\Bbbk)\cong\Bbbk\).
Hence \(\Gamma\) cannot be Cohen--Macaulay of dimension \(d\), ruling out
uniform Cohen--Macaulayness.

For double Cohen--Macaulayness, let \(q\) be the vertex of \(\Delta\)
corresponding to a maximal cube \(Q\), and let \(D\) be the subcomplex
of faces not containing \(q\). The link of \(q\) is
\(\operatorname{lk}_\Delta(q)=\operatorname{sd}(\partial Q)\),
a \((d-1)\)-sphere. Thus
\[
    H_d(\Delta,D;\Bbbk)
    \cong\widetilde H_{d-1}(\operatorname{lk}_\Delta(q);\Bbbk)
    \cong\Bbbk.
\]
Contractibility again gives \(\widetilde H_{d-1}(D;\Bbbk)\cong\Bbbk\),
so \(D\) cannot be Cohen--Macaulay of dimension \(d\).
This rules out double Cohen--Macaulayness as well.
\end{rem}
}

{
The preceding results concern the topology of the moment complex.
They do not by themselves determine the depth of the local rings of
\(\operatorname{QFl}_n\). Relating these conclusions to geometric
Cohen--Macaulayness requires information about the local equations and
scheme-theoretic intersections of its toric components.

\begin{ques}[Geometric Cohen--Macaulayness]\label{ques:geometric-cm}
Is the quasisymmetric flag variety \(\operatorname{QFl}_n\) locally Cohen--Macaulay for every \(n\ge1\)? Equivalently, is the local ring \(\mathcal O_{\operatorname{QFl}_n,x}\) Cohen--Macaulay at every point \(x\in\operatorname{QFl}_n\)?
\end{ques}
}

\section*{Acknowledgements}
We thank Hunter Spink, Patricia Klein, Vasu Tewari, Lucas Gagnon, and Philippe Nadeau for helpful conversations and suggestions.

\appendix
\section{Comparison with the moment complex}\label{app:comparison}

{
An indexed forest is a sequence of rooted plane full binary trees, all but
finitely many of them single leaves, with the leaves numbered consecutively.
Every internal node has exactly two children.
We use forests for which the leaves beyond \(n\) are single-leaf trees.
A top tree has all of \([n]\) as the leaves of its first tree.
Bergeron, Gagnon, Nadeau, Spink, and Tewari give a bijection
\(\operatorname{ForToNC}\) from these indexed forests on \([n]\) to
\(\operatorname{NC}_n\) \cite[Definition~7.10]{BergeronGagnonNadeauSpinkTewari2025}.
We use the characterization of its inverse by the forest code recalled below.

\begin{proof}[Proof of Proposition~\ref{prop:partition-model}]
Fix \(P\in\operatorname{NC}_n\). We first identify the coordinates of
\(Q_P\) with the internal nodes of the indexed forest corresponding to
\(P\). For an internal node \(v\) with left subtree \(L\), let
\(\tau_v=(i,j)\), where \(i\) is the rightmost leaf of \(L\) and \(j\)
is the rightmost leaf below \(v\). We construct the forest so that these
pairs are precisely the arcs and connectors in \(E(P)\).

Work from the innermost blocks outward. For a block
\(B=\{b_1<\cdots<b_m\}\), start with the leaf \(b_1\). Suppose the
tree constructed so far ends at \(b_r\). List the immediate children
\(C_1,\ldots,C_k\) of \(B\) in the gap \((b_r,b_{r+1})\) so that
\(\max C_h<\min C_{h+1}\) for \(1\le h<k\). Their trees occupy disjoint intervals and together fill the gap.
Start a new tree \(R\) with the leaf \(b_{r+1}\). For
\(h=k,k-1,\ldots,1\), replace \(R\) by the tree with a new root,
the tree of \(C_h\) as its left subtree, and the previous \(R\) as its
right subtree. The pair attached to this new root is
\((\max C_h,b_{r+1})\), the connector of \(C_h\). Now join the tree
ending at \(b_r\) and \(R\) as the left and right subtrees of a new
root. Its pair is \((b_r,b_{r+1})\), the next arc of \(B\).
Repeating this step gives a tree on all leaves from \(\min B\) to
\(\max B\), including the case of a singleton block. The trees of the
outer blocks, in order, followed by the single-leaf trees beyond \(n\),
form an indexed forest \(F\).

To identify \(F\), recall that its code entry \(c_i(F)\) counts the
internal nodes whose leftmost descendant leaf is \(i\). The \(|B|-1\)
nodes introducing arcs of \(B\) all have leftmost leaf \(\min B\).
For a nested block \(B\), the node introducing its connector contributes
one more at \(\min B\). These are all the nodes. Thus at \(i=\min B\)
the code entry is \(|B|\) for a nested block and \(|B|-1\) for an outer
block, with all other entries zero. This is the inverse-code description
of \(F_P=\operatorname{ForToNC}^{-1}(P)\) given immediately after
\cite[Definition~7.10]{BergeronGagnonNadeauSpinkTewari2025}.
The forest code is injective
\cite[Section~3.1]{BergeronGagnonNadeauSpinkTewari2025}, so \(F=F_P\).
In particular, \(v\mapsto\tau_v\) identifies its internal nodes with
\(E(P)\).

In the moment cube of \(F_P\), viewed as a forest with all internal
nodes black, the two opposite facets indexed by \(v\) correspond to
changing \(v\) to white and performing left-edge deletion at \(v\)
\cite[Proposition~5.4, Definition~5.6, and Theorem~5.7]{BergeronGagnonNadeauSpinkTewari2026}.
Assign these the coordinates \(x_{\tau_v}=0\) and \(x_{\tau_v}=1\),
respectively. At a vertex let \(S\) be the set of edges with
\(x_e=0\). By \cite[Theorem~6.3]{BergeronGagnonNadeauSpinkTewari2026},
the vertex permutation is the product of the transpositions in \(S\),
with ancestors preceding descendants. We compose permutations from right to left. Incomparable nodes give disjoint
transpositions, so their relative order does not affect the product.
Since \(G(P)\) is a forest by Lemma~\ref{lem:graph-forest}, each edge
joins two different components of the previously selected edges.
Inductively, multiplication by its transposition merges exactly the
two corresponding cycles. The final cycle partition is therefore
\(\pi_P(S)\). The permutation is noncrossing
\cite[Theorem~6.1]{BergeronGagnonNadeauSpinkTewari2026}, so its cycle
partition determines the permutation under our decreasing-cycle convention.
This identifies all the vertex labels of each coordinate face in \(Q_P\)
with those of the corresponding face of the moment cube.

Each indexed forest supported on \([n]\) occurs as \(F_P\) for a unique
\(P\). In particular the moment cubes of all the top trees occur, and
their faces exhaust \(\operatorname{Complex}(\operatorname{QFl}_n)\). Conversely, start with the first
tree of \(F_P\) on \([n]\) and join each successive tree on its right
under a new root. This gives a top tree. Left-edge deletion at the new
roots, starting with the last one added, separates their two subtrees
and recovers \(F_P\). By the facet description above, its moment cube
and all its coordinate faces therefore occur in the geometric complex.
By \cite[Theorem~8.1]{BergeronGagnonNadeauSpinkTewari2026}, one
geometric face is contained in another exactly when its vertex set is
contained in the other vertex set. Hence coordinate faces with equal
families of partition labels represent the same geometric face, and
inclusion of these families gives all incidences. This proves the
claimed isomorphism of face posets.

Finally, the all-black forest \(F_P\) is in normal form by
\cite[Definition~8.5]{BergeronGagnonNadeauSpinkTewari2026}. Hence
\cite[Proposition~8.10]{BergeronGagnonNadeauSpinkTewari2026}
gives \(\operatorname{ForToNC}(F_P)=P\) as the Bruhat maximum of the
moment cube of \(F_P\). For an arbitrary nonempty geometric face \(H\),
let \(\widehat H\) be its normal-form representative
\cite[Proposition~8.14]{BergeronGagnonNadeauSpinkTewari2026}.
Proposition~8.10 identifies the maximum of its vertex labels as the
partition assigned to \(\widehat H\) by the extension of
\(\operatorname{ForToNC}\) to normal forms
\cite[Section~8.2]{BergeronGagnonNadeauSpinkTewari2026}. Denote this maximum
by \(\lambda(H)\). Applying \cite[Lemma~9.3]{BergeronGagnonNadeauSpinkTewari2026} with
\(F=F_{\lambda(H)}\) shows that \(\widehat H\) is a face of this forest
cube. Under the isomorphism just proved, this gives
\(H\subseteq Q_{\lambda(H)}\).
\end{proof}

\begin{exam}
For \(P=13|2|4\), the construction gives the forest \(F_P\) on the left.
The graph \(G(P)\) is shown on the right.
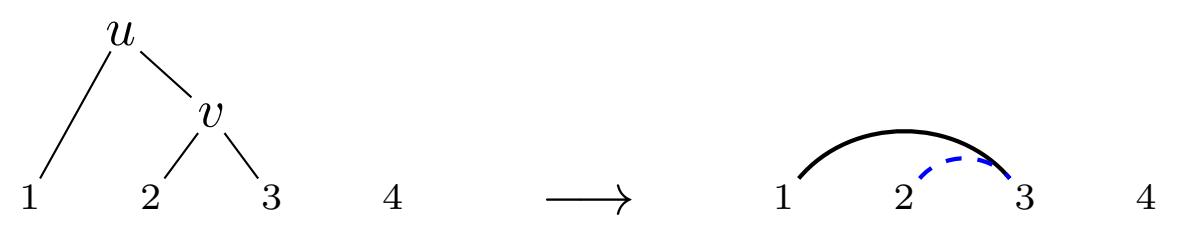
\[
\begin{tikzpicture}[x=0.8cm,y=0.6cm,color=black,
    every node/.style={inner sep=1pt,font=\small}]
    \foreach \i in {1,...,4}
        \node (l\i) at (\i,0) {$\scriptstyle \i$};
    \node (v) at (2.5,0.9) {$v$};
    \node (u) at (1.75,1.8) {$u$};
    \draw (l1) -- (u) -- (v) -- (l2) (v) -- (l3);
\end{tikzpicture}
\qquad\longrightarrow\qquad
\begin{tikzpicture}[x=0.8cm,y=0.8cm,color=black,
    every node/.style={inner sep=1pt,font=\small}]
    \foreach \i in {1,...,4}
        \node (l\i) at (\i,0) {$\scriptstyle \i$};
    \draw[thick] (l1) to[out=50,in=130] (l3);
    \draw[blue,dashed,thick] (l2) to[out=50,in=130] (l3);
\end{tikzpicture}
\]
The root \(u\) has left child \(1\) and rightmost leaf \(3\), so
\(\tau_u=(1,3)\), the arc of the block \(13\). The children of \(v\)
are \(2\) and \(3\), so \(\tau_v=(2,3)\), the connector of the
singleton \(2\). Joining \(u\) and leaf \(4\) under a new root gives
the top tree for \(\kappa(P)=134|2\) and adds the edge \((3,4)\).
\end{exam}
}

\bibliographystyle{amsrn}
\end{document}